\documentclass[12pt]{amsart}
\usepackage{amssymb,eepic,verbatim}
\newtheorem{theorem}{Theorem}[section]

\newtheorem{lemma}[theorem]{Lemma}
\theoremstyle{definition}

\theoremstyle{remark}
\newtheorem{remark}[theorem]{Remark}

\newcommand{\Z}{\mathbb{Z}}

\newcommand{\on}{\operatorname}

\newcommand{\Hom}{ \on{Hom}}

\newcommand{\ssm}{\kern-.5ex \smallsetminus \kern-.5ex}

\newcommand\dirac{/\kern-1.2ex\partial} 
\newcommand\qu{/\kern-.7ex/} 
 
\renewcommand{\comment}[1]   {{}}
\newcommand{\labell}\label

\newcommand{\ti}{\tilde}

\begin{document}

\def\mathunderaccent#1{\let\theaccent#1\mathpalette\putaccentunder}
\def\putaccentunder#1#2{\oalign{$#1#2$\crcr\hidewidth \vbox
to.2ex{\hbox{$#1\theaccent{}$}\vss}\hidewidth}}
\def\ttilde#1{\tilde{\tilde{#1}}}

\title{On the quantum product of Schubert classes}

\date{November 12, 2001}

\author{W. Fulton}\thanks{Partially supported by NSF grant 
DMS9970435}
\address{Department of Mathematics,
University of Michigan,
2074 East Hall,
Ann Arbor, Michigan  48109-1109, U.S.A.}
\email{wfulton@math.lsa.umich.edu}

\author{C. Woodward}\thanks{Partially supported by NSF grant 
DMS9971357}
\address{Mathematics-Hill Center, Rutgers University,
110 Frelinghuysen Road, Piscataway NJ 08854-8019, U.S.A.}
\email{ctw@math.rutgers.edu}

\begin{abstract}  We give a formula for the smallest powers of the quantum 
parameters  $q$  that occur in a product of Schubert classes in the (small) 
quantum cohomology of general flag varieties  $G/P$.  We also include a 
complete proof of Peterson's quantum version of Chevalley's formula, also for 
general  $G/P$'s.
\end{abstract}

\maketitle

\section{Introduction}

The Grassmannian was the first variety whose quantum cohomology was 
studied by physicists \cite{wi:vl}, and the first whose structure was worked 
out rigorously by mathematicians \cite{si:qc}, \cite{be:gi}, \cite{be:qs}.  
Other homogeneous varieties  $G/P$  have been studied (see below), but the 
story here remains far from complete.  Quantum cohomology has gone far 
beyond these beginnings, with all smooth projective varieties (or compact 
symplectic manifolds) enjoying a version of quantum cohomology.  However, 
there are still interesting questions to be answered about the case of  $G/P$  
in general, and Grassmannians in particular.  Our aim in this paper is to give 
an explicit formula for lowest degrees that occur in quantum product of 
Schubert classes.

The classical cohomology of a Grassmannian $\operatorname{Gr}(k,n)$ of
$k$-planes in $\mathbb{C}^n$ has a basis of Schubert classes
$\sigma_{\lambda}$, as $\lambda$ varies over partitions whose Young
diagram fits in a $k$ by $n-k$ rectangle.  The (complex) codimension
of $\sigma_{\lambda}$ is $|\lambda| = \sum \lambda_i$, the number of
boxes in the Young diagram.  The Littlewood-Richardson rule gives the
coefficients of a Schubert class $\sigma_{\nu}$ in a product
$\sigma_{\lambda}\cdot\sigma_{\mu}$, for $|{\nu}| = |\lambda| +
|{\mu}|$.  It is an easy and well-known fact that the classical
product $\sigma_{\lambda}\cdot\sigma_{\mu}$ is nonzero precisely when
$\lambda$ and the $180^{\circ}$ rotation of $ {\mu}$ fit in the $k$ by
$n-k$ rectangle without overlap; for example, the dual class to
$\sigma_{\lambda}$ is the class $\sigma_{\mu}$, for ${\mu} =
\lambda^{\vee}$ the partition such that $\lambda$ and the rotated
${\mu}$ exactly fill the rectangle without overlap.

The quantum cohomology of the Grassmannian is a free module over the 
polynomial ring  ${\mathbb{Z}}[q]$,  with a basis of Schubert classes; the 
variable $q$  has (complex) degree  $n$.  The quantum product  
$\sigma_{\lambda}\star\sigma_{\mu}$  is a finite sum of terms  $q^d 
\sigma_{\nu}$,  the sum over  $d\geq 0$  and  $|{\nu}| = |\lambda| + 
|{\mu}| - d\,n$,  each occurring with a nonnegative coefficient (a 
Gromov-Witten invariant); those with  $d = 0$  are the classical 
Littlewood-Richardson
coefficients.  This ring was studied in \cite{be:qm}, where an 
algorithm involving removing rim hooks was given for calculating these 
products.  It remains an important open problem to give a combinatorial 
formula for these coefficients (one that shows them to be nonnegative) when  
$d > 0$.

A simple argument due to Agnihotri showed that the quantum product  
$\sigma_{\lambda}\star\sigma_{\mu}$  of two Schubert classes in a 
Grassmannian can never be zero (see \cite{be:qm},\S5), so some  
$q^d\sigma_{\nu}$  must appear in such a product with positive coefficient.  
The problem we address here is to find the smallest power of  $q$  that 
occurs in a product  $\sigma_{\lambda}\star\sigma_{\mu}$.  The evidence 
from small examples, together with the role that rim hooks play in the 
quantum multiplication, lead one to conjecture that this smallest power of  
$d$  is the number of rim hooks it takes to cover the overlap of  $\lambda$  
and the $180^{\circ}$ rotation of  ${\mu}$  in the  $k$  by  $n-k$  
rectangle.  Equivalently,  $d$  is the maximum for which there is a diagonal 
sequence of boxes, from northwest to southeast, in this overlap.  Here is an 
example, for  $k = 4$,  $n = 9$,  $\lambda = (5, 4, 4, 3)$,  and  ${\mu} = 
(5,4,4,1)$:
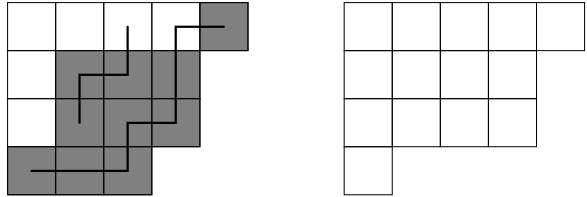
\begin{figure}[h]  
\setlength{\unitlength}{0.00083333in}
\begingroup\makeatletter\ifx\SetFigFont\undefined%
\gdef\SetFigFont#1#2#3#4#5{%
  \reset@font\fontsize{#1}{#2pt}%
  \fontfamily{#3}\fontseries{#4}\fontshape{#5}%
  \selectfont}%
\fi\endgroup%
{
\begin{picture}(3624,1239)(0,-10)
\path(12,1212)(312,1212)(312,912)
	(12,912)(12,1212)
\path(312,1212)(612,1212)(612,912)
	(312,912)(312,1212)
\path(612,1212)(912,1212)(912,912)
	(612,912)(612,1212)
\path(912,1212)(1212,1212)(1212,912)
	(912,912)(912,1212)
\path(12,912)(312,912)(312,612)
	(12,612)(12,912)
\path(12,612)(312,612)(312,312)
	(12,312)(12,612)
\path(2112,1212)(2412,1212)(2412,912)
	(2112,912)(2112,1212)
\path(2412,1212)(2712,1212)(2712,912)
	(2412,912)(2412,1212)
\path(2112,912)(2412,912)(2412,612)
	(2112,612)(2112,912)
\path(2112,612)(2412,612)(2412,312)
	(2112,312)(2112,612)
\texture{44555555 55aaaaaa aa555555 55aaaaaa aa555555 55aaaaaa aa555555 55aaaaaa 
	aa555555 55aaaaaa aa555555 55aaaaaa aa555555 55aaaaaa aa555555 55aaaaaa 
	aa555555 55aaaaaa aa555555 55aaaaaa aa555555 55aaaaaa aa555555 55aaaaaa 
	aa555555 55aaaaaa aa555555 55aaaaaa aa555555 55aaaaaa aa555555 55aaaaaa }
\shade\path(1212,1212)(1512,1212)(1512,912)
	(1212,912)(1212,1212)
\path(1212,1212)(1512,1212)(1512,912)
	(1212,912)(1212,1212)
\shade\path(312,912)(612,912)(612,612)
	(312,612)(312,912)
\path(312,912)(612,912)(612,612)
	(312,612)(312,912)
\shade\path(612,912)(912,912)(912,612)
	(612,612)(612,912)
\path(612,912)(912,912)(912,612)
	(612,612)(612,912)
\shade\path(912,912)(1212,912)(1212,612)
	(912,612)(912,912)
\path(912,912)(1212,912)(1212,612)
	(912,612)(912,912)
\shade\path(312,612)(612,612)(612,312)
	(312,312)(312,612)
\path(312,612)(612,612)(612,312)
	(312,312)(312,612)
\shade\path(912,612)(612,612)(612,612)
	(912,612)(912,612)
\path(912,612)(612,612)(612,612)
	(912,612)(912,612)
\shade\path(612,612)(912,612)(912,312)
	(612,312)(612,612)
\path(612,612)(912,612)(912,312)
	(612,312)(612,612)
\shade\path(912,612)(1212,612)(1212,312)
	(912,312)(912,612)
\path(912,612)(1212,612)(1212,312)
	(912,312)(912,612)
\shade\path(12,312)(312,312)(312,12)
	(12,12)(12,312)
\path(12,312)(312,312)(312,12)
	(12,12)(12,312)
\shade\path(312,312)(612,312)(612,12)
	(312,12)(312,312)
\path(312,312)(612,312)(612,12)
	(312,12)(312,312)
\shade\path(612,312)(912,312)(912,12)
	(612,12)(612,312)
\path(612,312)(912,312)(912,12)
	(612,12)(612,312)
\path(2712,1212)(3012,1212)(3012,912)
	(2712,912)(2712,1212)
\path(3012,1212)(3012,1212)(3012,1212)
	(3012,1212)(3012,1212)
\path(3012,1212)(3312,1212)(3312,912)
	(3012,912)(3012,1212)
\path(3312,1212)(3612,1212)(3612,912)
	(3312,912)(3312,1212)
\path(2412,912)(2712,912)(2712,612)
	(2412,612)(2412,912)
\path(2712,912)(3012,912)(3012,612)
	(2712,612)(2712,912)
\path(3012,912)(3312,912)(3312,612)
	(3012,612)(3012,912)
\path(2412,612)(2712,612)(2712,312)
	(2412,312)(2412,612)
\path(2712,612)(3012,612)(3012,312)
	(2712,312)(2712,612)
\path(3012,612)(3312,612)(3312,312)
	(3012,312)(3012,612)
\path(2112,312)(2412,312)(2412,12)
	(2112,12)(2112,312)
\thicklines
\path(1362,1062)(1062,1062)(1062,462)
	(762,462)(762,162)(162,162)
\path(762,1062)(762,762)(462,762)(462,462)
\end{picture}
}
\end{figure}

\noindent The overlap of $\lambda$ with the rotation of ${\mu}$ is
shaded, and one of the ways of covering the overlap with two rim hooks
is indicated.  In fact,\footnote{We recommend the program of Anders
Buch (http://home.imf.au.dk/abuch/lrcalc/) for computing classical and
quantum Littlewood-Richardson coefficients.}
$$
 \sigma_{5\,4\,4\,3}\star\sigma_{5\,4\,4\,1} = q^2(\sigma_{5\,3\,2\,2} + 
\sigma_{5\,3\,3\,1} + \sigma_{5\,4\,2\,1}) + q^3(\sigma_3 + 
2\sigma_{2\,1} + \sigma_{1\,1\,1}).  
$$
	
This conjecture is proved in this paper.\footnote{Buch \cite{bu:qc},
P. Belkale \cite{be:tr} and A. Yong \cite{yo:bo} have recently given
proofs of stronger versions of this result for the Grassmannians.}  In
fact, we prove a generalization of this conjecture for any $G/P$.  In
general the degree $d$ is an sequence of nonnegative integers, one for
each $1$-dimensional Schubert class.  We give a formula, in terms of
the combinatorics of the Bruhat order, for the smallest degrees $d$
such that $q^d$ occurs in a product of Schubert classes.  It follows
in particular that these quantum products can never be zero.

In the preprint version of this paper we posed the problem of giving a
criterion for exactly which powers of $q$ appear, or even for which
coefficients appear; for applications of such criteria, see
\cite{ag:ei}.  In the classical case, our understanding of this has
increased dramatically recently, thanks to Klyachko, Knutson, and Tao,
see \cite{fu:ev}.  For quantum products in Grassmannians, an upper
bound was obtained by A. Yong \cite{yo:bo}, and then the question of
which powers of $q$ appear was completely solved by A. Postnikov
\cite{po:af}.

An understanding of (small) quantum cohomology for a $G/P$ requires
first a presentation of the ring $QH^{\star}(G/P)$, and second, a
``quantum Giambelli formula'' for the class of a Schubert variety in
terms of this presentation.  This has been worked out for the variety
of complete flags (\cite{gi:qc}, \cite{cf:qc2}, \cite{fgp:qs},
\cite{ch:qf}) and partial flags (\cite{as:qc}, \cite{ki:pf},
\cite{cf:qc3}, \cite{ch:qf}), and recently for the Lagrangian
Grassmannian \cite{kt:lg}.  Descriptions of the quantum cohomology
ring have been given for general complete flag varieties $G/B$
\cite{ki:qc}, and partial descriptions for general $G/P$ by Peterson
\cite{pe:mi}, but Giambelli formulas are not yet known in general.

Most of what is known about quantum Giambelli formulas comes from 
computing formulas for degeneracy loci on certain Quot schemes (although 
A. Buch \cite{bu:qc} has recently given a proof for the Grassmannian, and 
Buch, Kresch, and Tamvakis \cite{bkt:lg} for some others, that does not 
depend on moduli spaces).  

Given this limited knowledge about quantum cohomology of general
$G/P$'s, it is somewhat surprising that we are able to solve this
problem for other $G/P$'s.  On the other hand, it indicates that, even
in type $A$, we will not use algebraic formulas for quantum Schubert
classes, and we will not use Quot schemes.  Rather, we use the spaces
$\overline{M}_{0,n}(X,d)$ of stable maps from genus $0$ curves with
$n$ marked points to $X = G/P$, which were constructed by Kontsevich
and Manin to prove the associativity of quantum products \cite{km:qh},
see \cite{fp:st}.

In the next few sections we lay out the necessary notation, recalling
the standard facts about the geometry of $G/P$'s and quantum
cohomology that are needed to state the theorem precisely.  Recall
that the Schubert classes $\sigma_u$ are parametrized by elements $u$
in $W/W_P$, where $W$ and $W_P$ are the Weyl groups of $G$ and $P$.
The idea behind one implication of this formula can be explained
roughly as follows.  If a product $\sigma_u\star\sigma_v$ contains a
term $q^d\sigma_w$, in the space $\overline{M}_{0,3}(X,d)$ the locus
of stable maps of degree $d$ to $X$ that map the first marked point to
a Schubert variety for $\sigma_u$ and the second marked point to an
opposite Schubert variety for $\sigma_v$ must contain a point that is
fixed by the maximal torus $T$ of $G$.  This fixed point is a map from
a curve $C$ into $X$, where $C$ is a tree of ${\mathbb{P}}^1$'s.  The
images of the intersections of the components of $C$ are fixed points
of $T$ in $X$, which are also indexed by elements of $W/W_P$.  This
produces a chain of elements of $W/W_P$, and this chain forces the
elements $u$ and $v$ to be close to each other in a certain way.  In
the case of the Grassmannian, this closeness translates to the
condition that the overlap described above can be covered by $d$ rim
hooks.

For the converse, any such chain does arise from a fixed point in such
a moduli space, but it is not obvious when a point corresponds to a
non-vanishing Gromov-Witten invariant.  The key to this is provided by
our transversality result in \S7, which we deduce from Kleiman's
general transversality theorem \cite{kl:gt}.  Similar ideas can be
used to prove Peterson's quantum extension of Chevalley formula for
multiplying a general Schubert class by a codimension one Schubert
class.  We have taken this opportunity to include a complete proof of
this formula in \S10.

The use of torus action in this setting goes back to Kontsevich
\cite{ko:lo}, who was inspired by Ellingsrud and Str{\o}mme.  It has
been used many times since, see \cite{ki:co}, \cite{th:ir}.  The idea
that Schubert varieties of opposite Borel subgroups are in general
position for the purposes of quantum cohomology we learned from
Peterson.

We thank Anders Buch and Alex Yong for several helpful conversations.
 
\section{Localization}

The following lemma, which is a special case of a theorem of Bott 
\cite{bo:rf}, provides simple proofs of the basic facts we need about divisors 
and curves on homogeneous varieties.

\begin{lemma} [Localization]  \label{lem1} Suppose a torus  $T$  acts on a 
curve  
$C \cong  {\mathbb{P}}^1$,  with fixed points  $p \neq   q$,  and suppose  
$L$  is a $T$-equivariant line bundle on  $C$.  Let  $\chi_ p$  and  
$\chi_q$  be the weights of  $T$  acting on the fibers  $L_p$  and  $L_q$,  
and let $\psi_  p$  be the weight of  $T$  acting on the tangent space to  
$C$  at $p$.  Then  
$$
\chi_p  -  \chi_q  =  n \, \psi_  p ,
$$
where  $n = \int_C  c_1(L)$  is the degree of  $L$.  
\end{lemma}

Note that $\psi_q = - \psi_p$,  so the result is independent of ordering of  
$p$  and  $q$.  Note also that both sides vanish if  $T$  acts trivially on  
$C$.

\section{Schubert varieties in  $G/P$}

We recall some basic notions about Schubert varieties and Schubert classes 
for a variety  $X = G/P$,  in order to fix our notation.  As usual,  $G$  
denotes a connected, simply connected, semisimple complex Lie group, in 
which we have fixed a Borel subgroup  $B$  and a maximal torus  $T$   in  
$B$.  We use the notation  $W$  for the Weyl group  $N(T)/T$,  $R = 
R^{+}\cup R^{-}$  for the roots (positive and negative), and  $\Delta$  for 
the simple roots; the {\bf reflections}  $s_{\alpha}$  in  $W$  are indexed by 
the positive roots  $\alpha$;  they are {\bf simple reflections} if  $\alpha$  
is in  $\Delta$.  The {\bf length}  $\ell(w)$  of an element  $w$  of  $W$  is 
the minimum number of simple reflections whose product is  $w$.  The 
element of longest length is denoted  $w_o$.  The {\bf opposite Borel} 
subgroup is $B\overline{\phantom{a}} = w_oBw_o$.

The parabolic subgroups  $P$  of  $G$  correspond canonically to subsets  
$\Delta_P$  of  $\Delta$.  Let  $R_P^{+}$  be the set of positive roots that 
can be written as sums of roots in  $\Delta_P$.  If  $\mathfrak{g} = 
\mathfrak{t}\, \oplus\, \bigoplus_{\alpha \in  R}  \mathfrak{g}_{\alpha}$  is 
the root space decomposition of the Lie algebra of  $G$,  then the Lie 
algebra  $\mathfrak{p}$  of  $P$  is the direct sum of $\mathfrak{t}$   and 
all  $\mathfrak{g}_{\alpha}$  for  $\alpha$  in  $R^{+} \cup (- R_P^{+})$.  
The group  $W_P$,  generated by the reflections  $s_{\alpha}$,  for  
$\alpha$  in  $\Delta_P$,  is the Weyl group of a Levi subgroup of  $G$  
corresponding to  $P$;  in particular,  $R_P^{+}$   is the corresponding set 
of positive roots, which consists of those  $\alpha$  in  $R^{+}$  such that  
$s_{\alpha}$  is in $W_P$.    

For an element  $u$  in  $W/W_P$,  $\ell(u)$  denotes the minimum length 
of a representative in  $W$.  In fact, each  $u$  has a unique representative 
of minimum length;  each element  $w$  of  $W$  can be written uniquely as 
a product  $a\cdot b$,  with $a$ the element of minimal length in the coset 
of  $w$ and  $b$  in  $W_P$,  and with  $\ell(w) = \ell(a) + \ell(b)$.  (For 
these facts see \cite{hu:rg}, \S1.10.)  The Weyl group acts on the left on  
$W/W_P$.  For  $u$  in  $W/W_P$,  we write  $u^{\vee}$  in place of 
$w_o u$.

For  $u$  in  $W/W_P$,  we let $ X(u) = \overline{BuP/P}$ be the 
corresponding {\bf Schubert variety}.  (The  $u$  on the right of this 
equation should be replaced by a representative first in  $W$,  and then by a 
representative in  $N(T)$,  but, as the result is independent of these 
choices, we follow the common convention of omitting them.)  This is a 
subvariety of  $X = G/P$  of dimension  $\ell(u)$;  we denote its cohomology 
class  $[X(u)]$  by  $\sigma(u)$.  Similarly, we let  $Y(u) = 
\overline{B\overline{\phantom{a}}uP/P}$  
 be the {\bf opposite Schubert variety}; it is of 
codimension  $\ell(u)$,  and we denote its cohomology class by  
$\sigma_u$.  Since  $Y(u) = w_o X(u^{\vee})$,  and translations of 
subvarieties by elements of  $G$  have the same cohomology classes, we have  
$$
\sigma_u  =  [Y(u)]  =  \sigma(u^{\vee})  =  [X(u^{\vee})]    \qquad \text {in} 
\qquad H^{2\ell(u)}(X).  
$$
(Cohomology here will always be taken with integer coefficients.) These 
classes form an additive basis for  $H^{\star}(X)$.   For  $w$  in  $W$,  we 
sometimes write  $\sigma_w$  for the class (and  $X(w)$  and  $Y(w)$  for 
the Schubert varieties) corresponding to the coset  $wW_P$  containing  $w$. 

For any $u$ in $W/W_P$, we let $x(u) = uP/P$ be the corresponding
point in $X$.  These are the fixed points of the action of $T$ on $X$.
The varieties $X(u)$ and $Y(u)$ meet transversally at the point
$x(u)$, and the classes $\sigma_u$ and $\sigma_{u^{\vee}} = \sigma(u)$
are dual classes under the intersection pairing: $\int_X
\sigma_u\cdot\sigma_v$ is $1$ if $v = u^{\vee}$, and $0$ otherwise.

The Schubert classes of dimension one have the form
$\sigma(s_{\beta})$ as $\beta$ varies over $\Delta \ssm \Delta_P$.  By a {\bf degree} $d$ we mean a nonnegative integral
combination $d = \sum d_{\beta}\, \sigma(s_{\beta})$ of these classes;
a degree may be identified with a collection of nonnegative integers
$(d_{\beta})_{\beta \in \Delta \ssm \Delta_P}$. The degrees
are the classes of curves on $X$.  If $d$ and $d'$ are degrees, we
write $d \leq d'$ to mean that $d_{\beta} \leq d'_{\beta}$ for all
$\beta$.

For any positive root  $\alpha$,  write  $\alpha = \sum  n_{\alpha \beta} 
\,  \beta$  as the nonnegative sum of simple roots  $\beta$;  then define the 
{\bf degree}  $d(\alpha)$  of  $\alpha$  by 
$$
d(\alpha)  =  \sum_{\beta  \in \Delta  \ssm  \Delta_P} \, 
n_{\alpha \beta} \, \frac{(\beta,\beta)} {(\alpha,\alpha)} \, 
\sigma(s_{\beta}).
$$
Here as usual $(\,\, , \,\,)$ is a $W$-invariant inner product on the
real subspace of $\mathfrak{t}^{\star}$ spanned by $R$.  If
$h_{\alpha} = 2\alpha/(\alpha,\alpha)$, and $\omega_{\beta}$ is the
fundamental weight corresponding to $\beta$ (so that the $h_{\beta}$
and $\omega_{\beta}$ form dual bases, for $\beta$ in $\Delta$), then
$h_{\alpha}(\omega_{\beta}) = n_{\alpha
\beta}(\beta,\beta)/(\alpha,\alpha)$, so this definition is equivalent
to setting
$$
d(\alpha)  =  \sum_{  \beta  \in   \Delta   \ssm   \Delta_P} \,  
h_{\alpha}(\omega_{\beta}) \, \sigma(s_{\beta}).  
$$

\begin{lemma} \label{lem2} If  $w$  is in  $W_P$,  then  $d(w(\alpha)) = 
d(\alpha)$.  
\end{lemma}

\begin{proof}  It suffices to prove this for a generator  $w = s_{\gamma}$,  
for  $\gamma$  in  $\Delta_P$.  Since  $s_{\gamma}(\alpha) = \alpha - 
2(\alpha,\gamma)/(\gamma,\gamma) \, \gamma$,  the coefficients of all  
$\beta$  in the expansions of  $\alpha$  and  $w(\alpha)$  are the same for  
$\beta$ in  $\Delta_P$.  Noting that  $(w(\alpha),w(\alpha)) = 
(\alpha,\alpha)$  for  all $w$  and  $\alpha$,  the result follows.  
\end{proof}

For any positive root  $\alpha$  that is not in  $R_P^{+}$,  there is a unique 
$T$-invariant curve  $C_{\alpha}$  in  $X$  that contains the points  $x(1)$  
and  $x(s_{\alpha})$.  Indeed,  $C_{\alpha} = Z_{\alpha}\cdot P/P$,  where  
$Z_{\alpha}$  is the $3$-dimensional subgroup of  $G$  whose Lie algebra is  
$\mathfrak{g}_{\alpha}\oplus\mathfrak{g}_{-\alpha} 
\oplus[\mathfrak{g}_{\alpha},\mathfrak{g}_{-\alpha}]$.  To see that  
$C_{\alpha}$  is unique, by the Bruhat decomposition there is a 
neighborhood of  $x(1)$  that is $T$-equivariantly isomorphic to  
$\mathfrak{u}^{-} = \mathfrak{g}/\mathfrak{p}$  (see \cite{bo:lag}, \S14); 
the $T$-invariant curves in  $\mathfrak{u}^{-}$  correspond to weight 
spaces  $\mathfrak{g}_{-\alpha}$  for  $\alpha$  in  $R^{+} \ssm 
R_P^{+}$.  If  $\alpha$  is in  $\Delta \ssm \Delta_P$,  then 
$C_{\alpha} = X(s_{\alpha})$  is one of our basic Schubert varieties.  If  
$\alpha$  is in  $R_P^{+}$,  then  $Z_{\alpha}\cdot P/P$  is the point  
$x(1)$.  

If  $\lambda$  is a weight that vanishes on all  $\beta$  in  $\Delta_P$,  it 
determines a character on  $P$,  and so a line bundle  $L(\lambda) = G 
\times_P \mathbb{C}(\lambda)$  on  $G/P$. 

\begin{lemma}  \label{lem3}  $\int_{C_{\alpha}} \,\,  c_1(L(\lambda)) = 
h_{\alpha}(\lambda)$.
\end{lemma}
\begin{proof} Lemma \ref{lem1} gives  $\int_{C_{\alpha}} \, 
c_1(L(\lambda)) = (\lambda - s_{\alpha}(\lambda))/\alpha = 
h_{\alpha}(\lambda)$.
\end{proof}

Applying this to  $\lambda = \omega_{\beta}$  and  $\alpha = \beta$  in  
$\Delta$,  we deduce:

\begin{lemma}  \label{lem4}  For  $\beta$  in  $\Delta \ssm 
\Delta_P$,  $c_1(L(\omega_{\beta})) = \sigma_{s_{\beta}}$. 
\end{lemma}

\begin{lemma}  \label{lem5} The degree  $[C_{\alpha}]$  of  $C_{\alpha}$  
is  $d(\alpha)$.
\end{lemma}

\begin{proof}  This is proved in \cite{ch:su}, pp. 14--19.  It follows more 
easily from the preceding two lemmas, since  $\sigma_{s_{\beta}}\cdot 
[C_{\alpha}] = h_{\alpha}(\omega_{\beta})$  implies that  $[C_{\alpha}] = 
\sum h_{\alpha}(\omega_{\beta}) \,  \sigma(s_{\beta})$.
\end{proof}

\begin{lemma}  \label{lem6}  The degree of the first Chern class of  $X$  on  
$C_{\alpha}$  is  $n_{\alpha} = 4(\rho_{P},\alpha)/(\alpha,\alpha)$,  where  
$\rho_{P} = \frac{1}{2} \sum  \gamma$,  with the sum over the positive 
roots  $\gamma$  not in  $R_P^{+}$.  In particular,  
$$
 c_1(T_X)  \,\,  =  \,\, 4 \sum_{\beta  \in  \Delta \ssm \Delta_P} 
\, \frac{(\rho_{P},\beta)}{(\beta,\beta)} \,  \sigma_{s_{\beta}} \,\,   = \,\, 2  
\sum_{\beta  \in   \Delta \ssm \Delta_P} h_{\beta}(\rho_{P}) \, 
\sigma_{s_{\beta}}.
$$
\end{lemma}

\begin{proof}  This can also be proved by localization.  Note that the tangent 
space to  $X$  at  $x(1)$  is  $\mathfrak{g}/\mathfrak{p} = 
\bigoplus_{\alpha \in  R'} \mathfrak{g}_{-\alpha}$,  where the sum is over 
the set  $R' = R^{+} \ssm R_P^{+}$.  So  $T_{x(u)}X  = 
\bigoplus_{\alpha \in  u\cdot R'}   \mathfrak{g}_{-\alpha}$.  The weight of 
the tangent space to  $C_{\alpha}$  at  $x(1)$  is  
$-\alpha$.  By Lemma \ref{lem1},  
$$
\int_{C_{\alpha}} \,   c_1(T_X) = \bigl(\sum_{\gamma  \in    R' } \,     
s_{\alpha}(\gamma) - \gamma \bigr)/(-\alpha) =  \sum_{\gamma  \in    R'} 
\,  2 \frac{(\gamma,\alpha)}{(\alpha,\alpha)}. 
$$
\end{proof}

\section{Chains in the Bruhat graph}

We need a combinatorial notion corresponding to the notion of a 
$T$-invariant curve joining the points  $x(u)$  and  $x(v)$  in  $X = G/P$.  

\begin{lemma} \label{lem7}  Let  $u$  and  $v$  be unequal elements in  
$W/W_P$.  The following are equivalent:
\begin{enumerate}
\item[(i)] \label{71} There is a reflection  $s$  in  $W$  
such that  $v = s\cdot 
u$.
\item[(ii)]  \label{72} There are representatives  $\tilde{u}$  for  $u$  and 
$\tilde{v}$  for  $v$  in  $W$,  and a reflection  $t$  in  $W$  such that 
$\tilde{v} = \tilde{u}\cdot t$. 
\item[(iii)] \label{73}  For any representative  $\tilde{u}$  of  $u$  in  $W$,  
there is a reflection  $s$  (resp. a reflection  $t$)  such that  $s\cdot 
\tilde{u}$  (resp.  $\tilde{u}\cdot t$)  is a representative of $v$.
\end{enumerate}

The reflection  $s$  of 
\textup{(i)} is uniquely determined.  The reflection  $t$  
of \textup{(ii)} is determined up to conjugation by an element of  $W_P$.   
\end{lemma}

\begin{proof}  (i) holds when there are representatives  $\tilde{u}$  
for $u$  and $\tilde{v}$  for  $v$  such that $\tilde{v} = s\cdot \tilde{u}$.  
Equivalently $\tilde{v} = \tilde{u}\cdot t$,  with  $t = \tilde{u}^{-1}\cdot 
s\cdot \tilde{u}$,  which is (ii).  In either case the representative  
$\tilde{u}$  can be chosen arbitrarily.  Both uniqueness assertions will follow 
from the 

\vspace{12pt}

\noindent
{\bf Claim.}  Let  $s$  and  $s'$  be reflections, not in  $W_P$.  If  $s' = 
s\cdot a$  for  some  $a$  in  $W_P$,  then  $s' = s$. 

\vspace{12pt}

Granting the claim, in (i), if  $s'\cdot u = s\cdot u \neq   u$,  then for 
any representative  $\tilde{u}$  of  $u$,  $s'\cdot \tilde{u} = s\cdot 
\tilde{u}\cdot a$  for some  $a$  in  $W_P$.  Then  
$\tilde{u}^{-1}\cdot s'\cdot\tilde{u} = \tilde{u}^{-1}\cdot 
s\cdot\tilde{u}\cdot a$,  and since  $\tilde{u}^{-1}\cdot s\cdot\tilde{u}$  is 
not in  $W_P$,  the claim implies that  $\tilde{u}^{-1}\cdot s'\cdot\tilde{u} 
= \tilde{u}^{-1}\cdot s\cdot \tilde{u}$,  so $s' = s$.  Similarly in (ii), if 
$\tilde{v} = \tilde{u}\cdot t$  and $\tilde{v}\cdot a =  \tilde{u}\cdot b\cdot 
t'$  for some  $a$  and  $b$  in  $W_P$,  then  $b\cdot t' = t\cdot a$,  so  
$b\cdot t'\cdot b^{-1} = t\cdot c$,  with $c = a\cdot b^{-1}$  in  $W_P$.  
The claim implies that  $t = b\cdot t'\cdot b^{-1}$,  as required.  

\vspace{12pt}
\noindent
{\bf Proof of the claim.}  Let  $v$  be a weight such that  $(\beta,v) = 0$  
for 
all  $\beta$  in  $\Delta_P$,  and  $(\beta,v) > 0$  for all  $\beta$  in 
$\Delta \ssm \Delta_P$.  For any  $w$  in  $W$,  we have  $w(v) = 
v$  if and only if  $w$  is in  $W_P$  (\cite{bo:lg} V, \S4.6).  In particular  $s'(v) = 
s(v) \neq   v$.  If  $s = s_{\alpha}$  and  $s' = s_{\gamma}$,  for  $\alpha$  
and  $\gamma$ positive roots, then  $v - 
2(\alpha,v)/(\alpha,\alpha)\cdot\alpha = v - 
2(\gamma,v)/(\gamma,\gamma)\cdot\gamma$.  This implies that  $\alpha$  
and  $\gamma$  are proportional, which cannot happen unless  $\alpha  = 
\gamma$.  
\end{proof}

The claim amounts to the fact that if  $\alpha$  and  $\gamma$  are distinct 
positive roots that are not sums of roots in  $\Delta_P$,  then the cosets of  
$s_{\alpha}$  and $s_{\gamma}$  in  $W/W_P$  are distinct. 

We will say that two unequal elements  $u$  and  $v$  in  $W/W_P$  are {\bf 
adjacent} if they are related as in Lemma \ref{lem7}.  Note that this is a 
symmetric relation.  In this case we define  $d(u,v)$  to be the degree  
$d(\alpha)$,  where  $t = s_{\alpha}$  is a reflection relating them as in 
(ii).  If  $t$  is replaced by  $w\cdot t\cdot w^{-1} = s_{w(\alpha)}$,  
for  $w$  in  $W_P$,  the degree does not change (Lemma \ref{lem2}), so  
$d(u,v)$  depends only on  $u$  and  $v$.  Note that if  $u$  and  $v$  are 
adjacent, then for any  $w$  in  $W$,  $w\cdot u$  and  $w\cdot v$  are also 
adjacent, and  $d(w\cdot u, w\cdot v) = d(u,v)$.  In particular, $u^{\vee}$  
and  $v^{\vee}$  are also adjacent, with $ d(u^{\vee},v^{\vee}) = d(u,v)$.

\begin{lemma} \label{lem8}   Elements  $u$  and  $v$  in  $W/W_P$  are 
adjacent if and only if  $x(u) \neq   x(v)$  and there is a $T$-invariant curve  
$C$  containing  $x(u)$  and  $x(v)$.  If this is true, the curve  $C$  is 
unique, isomorphic to  ${\mathbb{P}}^1$,  and its class  $[C]$  in  $H_2(X)$  
is equal to $d(u,v)$.
\end{lemma}

\begin{proof}  We have seen that the $T$-invariant curves containing  
$x(1)$  are exactly the curves  $C_{\alpha}$,  which also contains  
$x(s_{\alpha})$,  for  $\alpha$  in  $R^{+} \ssm R_P^{+}$ .  
General $T$-invariant curves in  $X$  therefore have the form  $w\cdot 
C_{\alpha}$,  for some  $\alpha$  in  $R^{+} \ssm R_P^{+}$   and  
$w$  in  $W$.  This curve is the unique $T$-invariant curve containing 
$x(w) = w\cdot x(1)$  and  $x(w\cdot s_{\alpha}) = w\cdot x(s_{\alpha})$.  
The result then follows from Lemmas \ref{lem2}, \ref{lem5}, and 
\ref{lem7}.  (For more about $T$-invariant curves in general, see 
\cite{ca:de}.)  
\end{proof}

We use also the {\bf Bruhat order} on  $W/W_P$,  which sets  $u \prec v$  if  
$X(u) \subset X(v)$.  

\begin{lemma} \label{lem9}  For  $u$  and  $v$  in  $W/W_P$,  the 
following are equivalent:
\begin{enumerate}
\item[(i)] \label{91}  $u \preceq v$;  
\item[(ii)] \label{92}  for any sequence of  
$\ell(v)$  simple reflections whose product represents  $v$,  a 
representative of  $u$  can be obtained by removing some of these 
transpositions;  \item[(iii)] \label{93}  $x(u) \in  X(v)$;  
\item[(iv)] \label{94} $v^{\vee} \preceq u^{\vee}$;  \item[(v)] \label{95}  $x(v) \in  
Y(u)$.
\end{enumerate}
\end{lemma}

\begin{proof}  For the equivalence of (i), (ii), and (iv), see 
\cite{hu:rg}, \S5.9, 5.10.  The equivalence of (i) and (iii) follows 
from the fact that  $X(u)$  is the closure of  $B\cdot x(u)$.  Then (iv) is 
equivalent to  $x(v^{\vee})$  being in  $X(u^{\vee})$,  or to  $x(v) = w_o 
x(v^{\vee})$  being in  $w_o X(u^{\vee}) = Y(u)$,  which is (v).
\end{proof}

Now define a {\bf chain} from  $u$  to  $v$  in  $W/W_P$  to be a sequence  
$u_0, u_1, \ldots , u_r$  in  $W/W_P$  such that  $u_i$  and  $u_{i-1}$  are 
adjacent for  $1 \leq i \leq r$,  and, in addition,  $u \preceq u_0$  and  $u_r 
\preceq v^{\vee}$.  For any chain  $u_0, u_1, \ldots , u_r$  we define the {\bf 
degree} of the chain to be the sum of the degrees  $d(u_{i-1},u_i)$,  for  $1 
\leq i \leq r$.  Note that such a chain from  $u$  to  $v$  determines a chain 
from  $v$  to  $u$,  by  $v \preceq {u_r}^{\vee}, \ldots , {u_0}^{\vee} \preceq 
u^{\vee}$,  and these chains have the same degree.  
Note also that there is a chain of degree  $0$  between  $u$  and  $v$  exactly 
when  $u \preceq v^{\vee}$.   

A chain from $u$  to  $v$  determines, and is determined by, a sequence of 
$T$-invariant curves  $C_1, C_2, \ldots , C_r$  in  $X$,  each meeting the 
next, with  $C_1$  meeting  $Y(u)$  and  $C_r$  meeting  $X(v^{\vee})$.  
Indeed,  $C_i$  is the $T$-invariant curve that connects  $x(u_{i-1})$  to  
$x(u_i)$.  The degree of the chain is the sum of the classes  $[C_i]$  of the 
curves. 

\section{Interpretation for the Grassmannian}

For the Grassmannian  $\operatorname{Gr}(r,n)$,  $G$  is  
$SL_n(\mathbb{C})$,  $B$  is the subgroup of upper triangular matrices, 
$T$   the diagonal matrices in  $B$,  and  $W$  is identified with the 
symmetric group  $S_n$.  The simple roots are  $\Delta = \{\alpha_i = e_i - 
e_{i+1},  1 \leq i \leq n-1\}$.  The parabolic subgroup  $P$  consists of 
matrices in  $G$  that map the subspace of  $\mathbb{C}^n$  spanned by 
the first  $r$  basic vectors to itself, and  $\Delta_P$  consists of all simple 
roots with the exception of  $\alpha_r$.  The Weyl group  $W_P$  is 
identified with  $S_r \times S_{n-r}$.  The minimal representative of a  $u$  
in  $W/W_P$  is a permutation  $w$  such that  $w(1) < w(2) < \ldots < 
w(r)$  and  $w(r+1) < \ldots < w(n)$.  From this we form a partition
$$
\lambda(u)  =  (w(r) - r, w(r-1) - (r-1), \ldots , w(2) - 2, w(1) - 1),
$$
with  $n-r \geq \lambda_1(u) \geq \ldots \geq \lambda_r(u) \geq 0$.  This 
sets up natural bijections between:  (i) elements of  $W/W_P$;  (ii)  
partitions inside the  $r$  by  $n-r$  rectangle; and (iii) subsets of  $\{1, 
\ldots , n\}$  with  $r$  elements.  With this notation, the subset 
corresponding to  $u$  is  $\{w(1), w(2), \ldots , w(r)\}$,  for any 
representative  $w$  of  $u$  in  $S_n$.  Note that  $u \preceq v$  if and only 
if  $\lambda(u)$  is contained in  $\lambda(v)$,  i.e.,  $\lambda_i(u) \leq 
\lambda_i(v)$  for  $1 \leq i \leq r$.  
	
\begin{lemma} \label{lem10}  In the Grassmann case,  $u$  is adjacent to  
$v$  if and only if one of  $\lambda(u)$  and  $\lambda(v)$  is contained in 
the other, and the difference is a (connected, nonempty) rim hook.  In this 
case the degree  $d(u,v)$  is  $1$.
\end{lemma}

\begin{proof}  Let  $I$  and  $J$  be the subsets of  $\{1, \ldots , n\}$  
corresponding to $u$  and  $v$.  Then  $u$  and  $v$  are adjacent exactly 
when  $I$  and  $J$  differ by one element, i.e., there is a  $p$  in  $I 
\ssm J$  and a  $q$  in  $J \ssm I$  with  $I\cup{q} = 
J\cup{p}$;  in this case, $v = (p,q)\cdot u$.  If  $p < q$,  then  
$\lambda(u)$  is obtained from  $\lambda(v)$  by removing a rim hook of  
$q - p$  boxes, starting at the end the  $k^{\text {th}}$  row, where  $k - 1$  
is the number of elements in  $J$  that are bigger than  $q$.  The 
transposition  $t$  has the form  $(i,j)$  for $i \leq r < j$,  and one sees 
readily that  $d(e_i - e_j) = 1$.
\end{proof}

The Bruhat graph for this case is known as the {Johnson graph}.
The figure shows the Johnson graph for $\operatorname{Gr}(2,4)$,
labeled by the partitions $\lambda$ with $2\geq \lambda_1 \geq
\lambda_2 \geq 0$.
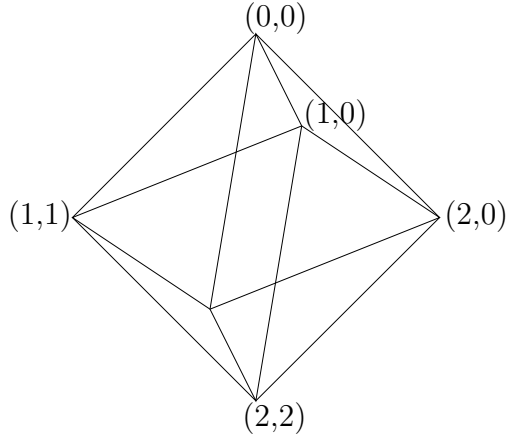
\begin{figure}[h]
\setlength{\unitlength}{0.00043333in}
\begingroup\makeatletter\ifx\SetFigFont\undefined%
\gdef\SetFigFont#1#2#3#4#5{%
  \reset@font\fontsize{#1}{#2pt}%
  \fontfamily{#3}\fontseries{#4}\fontshape{#5}%
  \selectfont}%
\fi\endgroup%
{
\begin{picture}(5334,4905)(0,-10)
\path(2640,4710)(390,2460)(2640,210)
	(4890,2460)(2640,4710)
\path(3203,3585)(4890,2460)(2078,1335)
	(390,2460)(3203,3585)
\path(2640,4710)(3203,3585)(2640,210)
	(2078,1335)(2640,4710)
\put(2490,4810){\makebox(0,0)[lb]{\smash{{{\SetFigFont{12}{14.4}{\rmdefault}{\mddefault}{\updefault}(0,0)}}}}}
\put(4950,2385){\makebox(0,0)[lb]{\smash{{{\SetFigFont{12}{14.4}{\rmdefault}{\mddefault}{\updefault}(2,0)}}}}}
\put(2475,-100){\makebox(0,0)[lb]{\smash{{{\SetFigFont{12}{14.4}{\rmdefault}{\mddefault}{\updefault}(2,2)}}}}}
\put(-400,2385){\makebox(0,0)[lb]{\smash{{{\SetFigFont{12}{14.4}{\rmdefault}{\mddefault}{\updefault}(1,1)}}}}}
\put(3225,3615){\makebox(0,0)[lb]{\smash{{{\SetFigFont{12}{14.4}{\rmdefault}{\mddefault}{\updefault}(1,0)}}}}}
\end{picture}
}
\end{figure}

\begin{remark}  In this Grassmann case, it is not hard to see that if there is 
a chain of degree  $d$  from  $u$  to  $v$,  then there is a {\em monotone} 
chain of degree at most  $d$  from  $u$  to  $v$,  i.e., a chain that starts 
with  $u_0$  equal to  $u$,  removes a rim hook at each stage, and ends at a  
$u_r$  with  $\lambda(u_r)$ contained in the $180^{\circ}$ rotation of  
$\lambda(v)$.  
\end{remark}

\section{The quantum cohomology of  $G/P$}

We next describe the (small) quantum cohomology of  $X$.  Take a variable  
$q_{\beta}$  for each  $\beta$  in  $\Delta \ssm \Delta_P$,  and 
let  ${\mathbb{Z}}[q]$  be the polynomial ring with these  $q_{\beta}$  as 
indeterminants, but giving  $q_{\beta}$  the degree  $2n_{\beta}$,  where  
$n_{\beta} = \int_{\sigma(s_{\beta})} c_1(T_X)$  from Lemma \ref{lem6}.  
For a degree $ d  =  \sum d_{\beta}\,\sigma(s_{\beta})$,  we write  $q^d$  
for the monomial  $\prod_{\beta}  q_{\beta}^{d_{\beta}}$.  The small 
quantum cohomology ring  $QH^{\star}(X)$  is, as a  
${\mathbb{Z}}[q]$-module,
simply  $H^{\star}(X)\otimes_{\mathbb{Z}}{\mathbb{Z}}[q]$,  so 
the same Schubert classes  $\sigma_u = \sigma_u {\small \otimes}  1$  
form a basis  for  $QH^{\star}(X)$  over  ${\mathbb{Z}}[q]$.  The 
multiplication is a deformation of the classical multiplication:
$$
\sigma_u  \star  \sigma_v  \, = \,  \sum_d  q^d \sum_w  \, N_{u,v}^w(d) \,  
\sigma_w, 
$$
where the first sum is over all degrees $d$, and the second is over
all $w$ in $W/W_P$ such that $\ell(w) = \ell(u) + \ell(v) -
\sum_{\beta} d_{\beta}n_{\beta}$.  The coefficient $N_{u,v}^w(d)$ is a
Gromov-Witten (GW) invariant: it is the number of morphisms $\varphi :
{\mathbb{P}}^1 \to X$ of degree $d$ (i.e.,
$\varphi_{\star}[{\mathbb{P}}^1] = d$ in $H_2(X)$), such that, for
three given distinct points $p_1, p_2, p_3$ in ${\mathbb{P}}^1$, and
three general $g_1, g_2$, and $g_3$ in $G$, $\varphi(p_1)$ is in
$g_1\cdot Y(u)$, $\varphi(p_2)$ is in $g_2\cdot Y(v)$, and
$\varphi(p_3)$ is in $ g_3\cdot X(w)$.  When $d = 0$, this is the
usual coefficient of $\sigma_w$ in the classical product
$\sigma_u\cdot\sigma_v$, which is the same as the intersection number
$\int_X \sigma_u\cdot\sigma_v\cdot\sigma_{w^{\vee}}$.

More generally, the (small) GW-invariant $\langle \sigma_{u_1},
\sigma_{u_2}, \ldots , \sigma_{u_n}\rangle_d$ can be defined whenever
$\sum \ell(u_i) = \operatorname{dim}(X) + \sum_{\beta} n_{\beta}
d_{\beta}$.  Fix general distinct points $p_1, \ldots , p_n$ in
${\mathbb{P}}^1$.  This invariant is the number of maps $\varphi$ from
${\mathbb{P}}^1$ to $X$ such that $\varphi(p_i)$ is in $ g_i\cdot
Y(u_i)$, for $1 \leq i \leq n$, and $g_1, \ldots , g_n$ general
elements of $G$.  This can be interpreted in the cohomology of
appropriate moduli spaces.  As we will need these spaces in our
proofs, we describe them now.

Let  $\overline{M}_{0,n}(X,d)$  be the moduli space of stable maps of 
degree  $d$  of   $n$-pointed genus  $0$  curves into  $X$;  a point is 
written  $(C, p_1, \ldots , p_n, \varphi)$,  where  $C$  is a connected tree 
of projective lines, meeting in nodes,  $p_1, \ldots , p_n$  are distinct 
nonsingular points of  $C$,  and  $\varphi : C \to  X$  is a morphism with  
$\varphi_{\star}[C] = d$, with the property that any component of  $C$  
that is mapped to a point by  $\varphi$  must have at least three points that 
are either marked points or intersection points with other components.  
This moduli space is a projective variety of dimension  
$$
\operatorname{dim}( \overline{M}_{0,n}(X,d)) \, =\,  \operatorname{dim}( X) 
\,+\, \sum_{\beta} n_{\beta} \, d_{\beta}\, + \,n \, -\, 3.  
$$
It comes equipped with  $n$  evaluation maps  $e_i  : 
\overline{M}_{0,n}(X,d) \to  X$,  taking  $(C, p_1, \ldots , p_n, \varphi)$  
to  $\varphi(p_i)$,  and a forgetful map  
$f  :  \overline{M}_{0,n}(X,d) \to  \overline{M}_{0,n}$,  where the latter is 
the space of stable $n$-pointed curves of genus $0$;  $f$  takes  $(C, p_1, 
\ldots , p_n, \varphi)$  to  $(C, p_1, \ldots , p_n)$,  but suitably stabilized 
by collapsing components of  $C$  that have fewer than three markings or 
intersections with other components.   We refer to \cite{fp:st} for 
construction and basic properties of these spaces and mappings, as well as 
Kontsevich's proof of the associativity of the quantum product.  

The GW-invariant  $\langle \sigma_{u_1}, \ldots , \sigma_{u_n}\rangle_d$   
is then the intersection number  $\int_{\overline{M}_{0,n}(X,d)} 
f^{\star}([p])\cdot {e_1}^{\star}(\sigma_{u_1})\cdot 
{e_2}^{\star}(\sigma_{u_2})\cdot \ldots \cdot 
{e_n}^{\star}(\sigma_{u_n})$,  where  $p$  is a point in  
$\overline{M}_{0,n}$.  Equivalently, it is the coefficient of the fundamental 
class  $1 = [\overline{M}_{0,n}]$  in the class
$$
f_{\star} \bigl( {e_1}^{\star}(\sigma_{u_1})\cdot 
{e_2}^{\star}(\sigma_{u_2})\cdot \ldots \cdot {e_n}^{\star}(\sigma_{u_n}) 
\bigr)   \qquad \text {in} \quad  H^0(\overline{M}_{0,n}).
$$
In particular, this shows that the Gromov-Witten invariants are the same 
whether one chooses any distinct points $p_1, \ldots , p_n$  in  
${\mathbb{P}}^1$  instead of general points (see \cite{be:qs}).  

The coefficient  $N_{u,v}^w(d)$  is equal to  $\langle \sigma_u, \sigma_v, 
\sigma_{w^{\vee}}\rangle_d$ .  In fact, these invariants can be used to 
multiply several Schubert classes directly:  
$$
\sigma_{u_1}\star \sigma_{u_2}\star \ldots \star \sigma_{u_n} =  \sum_d 
q^d \sum_w \, \langle \sigma_{u_1}, \ldots , \sigma_{u_n}, 
\sigma_{w^{\vee}} \rangle_d \, \sigma_w.
$$

\section{Transversality}

The results of this section are the main tools needed to prove our
theorem.  We give here a simple proof based on Kleiman's
transversality theorem.  An alternative proof is sketched briefly at
the end of this section.

\begin{lemma} \label{lem11}   Let  $U \subset G \times G$  be open, nonempty, 
and invariant under the left diagonal multiplication by $G$.  Let  $u_1$  and  
$u_2$ be in  $W/W_P$.  Then for any  $g_1, g_2$  in  $G$  such that 
$g_1B{g_1}^{-1}$ and $g_2B{g_2}^{-1}$  intersect in a maximal torus, 
there is a  $(h_1,h_2)$  in  $U$ such that
$$
h_1X(u_1) = g_1X(u_1)  \qquad \text  {and} \qquad    h_2X(u_2) = 
g_2X(u_2).
$$
\end{lemma}

\begin{proof}  Take  $(h_1,h_2)$ in the intersection of  $U$  with the open
set of pairs $(h_1,h_2)$ such that $h_1B{h_1}^{-1}$ and
$h_2B{h_2}^{-1}$ are opposite Borels, i..e., $h_1B{h_1}^{-1} \cap
h_2B{h_2}^{-1}$ is a maximal torus.  By \cite{bo:lag}, \S14.1, Cor. 3,
there is a $g$ in $G$ such that $g(h_iB{h_i}^{-1})g^{-1} =
g_iB{g_i}^{-1}$ for $i = 1, 2$.  Since $B$ is its own normalizer, this
implies that $gh_iB = g_iB$ for $i = 1, 2$.  Then $(gh_1,gh_2)$ is in
$U$, and $gh_iX(u_i) = g_iX(u_i)$ for $i = 1, 2$.
\end{proof}

\begin{lemma} \label{lem12}  Let  $Z$  be an irreducible  $G$-variety, and 
let  $F : Z \to  X \times X$  be a $G$-equivariant morphism, where  $G$  
acts diagonally on  $X \times X$.  Then, for any  $u$  and  $v$  in  
$W/W_P$,  the subscheme  $F^{-1}(Y(u) \times X(v))$  is reduced, locally 
irreducible, of codimension  $\ell(u) + \ell(v^{\vee})$,  and nonsingular at 
any nonsingular point of  $Z$  that maps to a nonsingular point of  $Y(u) 
\times X(v)$.  
\end{lemma}

\begin{proof}  Consider the diagram 
\[ \begin{array}{ccc}
& & X(u^{\vee}) \times X(v) \\
& & \downarrow \\
Z  &  \longrightarrow   &  X \times X
\end{array} \]
with  $G \times G$  acting on  $X \times X$.  Kleiman's transversality 
theorem \cite{kl:gt} produces a nonempty open set  $U$  of pairs  
$(g_1,g_2)$  in  $G \times G$  such that  $F^{-1}(g_1X(u^{\vee}) \times 
g_2X(v))$  satisfies the conclusions of the lemma.  (The dimension assertion 
would be valid in all characteristics; the others use the characteristic zero 
assumption.)  This set  $U$  is invariant by the left diagonal action of  $G$  
because the morphism  $F$  is $G$-equivariant. 

Now apply Lemma \ref{lem11} to  $(u_1,u_2) = (u^{\vee},v)$  and  
$(g_1,g_2) = (w_o,1)$.  This produces a  $(h_1,h_2)$  in  $U$  such that  
$h_1X(u^{\vee}) = w_oX(u^{\vee}) = Y(u)$  and  $h_2X(v) = X(v)$,  and 
Lemma \ref{lem12} follows.
\end{proof}

A point  $\zeta = (C, p_1, \ldots , p_n, \varphi)$  in a moduli space  
$\overline{M}_{0,n}(X,d)$  of stable maps consists of a tree  $C$ of  
${\mathbb{P}}^1$'s,  and marked points on some of its components, with a 
stable map from $C$ to  $X$.  Such a point lies in a unique locally closed 
subscheme  $V$,  for which the tree of curves has the same topological type 
(or combinatorial configuration), with marked points on corresponding 
components (see \cite{bm:st}, \cite{fp:st}).  The codimension of  $V$  is the 
number of nodes of  $C$.   

When  $n$  and  $d$  are understood, we set
$$
E(u,v)  =  {e_1}^{-1}(Y(u)) \cap {e_2}^{-1}(X(v)) ,
$$
a closed subscheme of  $\overline{M}_{0,n}(X,d)$. 

\begin{lemma} \label{lem13}  If  $E(u,v)$  is not empty, then  $E(u,v)$  is a 
reduced, locally irreducible, subscheme of  $\overline{M}_{0,n}(X,d)$,  of 
pure codimension  $\ell(u) + \ell(v^{\vee})$,  any component of which 
which meets any stratum  $V$  properly.  In particular, each irreducible 
component of  $E(u,v)$ meets the locus  $M_{0,n}(X,d)$  consisting of 
those  $(C, p_1, \ldots , p_n, \varphi)$  with  $C \cong {\mathbb{P}}^1$.
\end{lemma}

\begin{proof}  This follows from Lemma \ref{lem12}, and the fact that the 
strata are locally closed, $G$-invariant subvarieties in  
$\overline{M}_{0,n}(X,d)$.  
\end{proof}

\begin{lemma} \label{lem14}  For any degree  $d$,  for  $n > \sum  
n_{\beta}d_{\beta}$,  and for any distinct points  $p_1, \ldots , p_n$  in  
${\mathbb{P}}^1$  and any points  $x_1, \ldots , x_n$  in  $X$,  there are 
only finitely many morphisms  $\varphi : {\mathbb{P}}^1 \to  X$  of degree  
$d$  with  $\varphi(p_i) = x_i$  for  $1 \leq i \leq n$. 
\end{lemma}

\begin{proof}   Let $\operatorname{Hom}({\mathbb{P}}^1,X)_d$  
be the space of morphisms 
of degree  $d$  from  ${\mathbb{P}}^1$  to  $X$,  and let  
$e : \operatorname{Hom}({\mathbb{P}}^1,X)_d 
\to  X^n$ be the morphism obtained by 
evaluating at the given points  $p_1, \ldots , p_n$.  We show that $e$ is 
unramified, and hence has finite fibers.  This will be true if its tangent map
\cite{ko:ra}
$$
\Gamma({\mathbb{P}}^1,\varphi^{\star}(T_X))  \longrightarrow 
\oplus_{i = 1}^n \varphi^{\star}(T_X)(p_i) = \oplus_{i = 1}^n T_{x_i}(X)
$$
is injective.  Now  $\varphi^{\star}(T_X) = \oplus \mathcal{O}(m_j)$,  with  
$m_j \geq 0$,  and  $\sum m_j = \sum  n_{\beta}d_{\beta} < n$,  so  $m_j 
< n$  for all $j$.  The lemma follows from the elementary fact that  
$\Gamma({\mathbb{P}}^1,\mathcal{O}(m)) \to \oplus_{i = 1}^n  
\mathcal{O}(m)(p_i)$  is injective for  $m < n$  (since a nonzero polynomial 
of degree  $m$  cannot vanish at more than  $m$  points).
\end{proof}

The transversality Lemma \ref{lem13} can also be proved by showing that, 
given a fixed point  $\zeta = (C, p_1, \ldots , p_n, \varphi)$  with  
$\varphi(p_1) = x(u)$  and  $\varphi(p_2) = x(v)$,  the map from 
$\Hom(C,X)_d$  to  $X^2$ given by the first two evaluation maps is 
transversal to the subvariety $Y(u) \times X(v)$  at the point $x(u) \times 
x(v)$.  One shows that the map from the tangent space 
$\Gamma(C,\varphi^{\star}(T_X))$ to the normal space to $Y(u) \times 
X(v)$ at $x(u) \times x(v)$ is surjective.  This can be achieved by 
$T$-equivariantly decomposing the bundle $\varphi^{\star}(T_X)$ as a 
direct sum of line bundles, and using the fact that the weights of the normal 
spaces of  $Y(u)$ at $x(u)$ and $X(v)$ at $x(v)$ are disjoint.

\section{Chevalley's formula in the parabolic case}

Chevalley's formula \cite{ch:su} generalizes Monk's formula from the 
classical flag variety to an arbitrary  $G/B$,  giving a formula for the product 
of a codimension one Schubert class  $\sigma_{s_{\beta}}$  and an arbitrary 
Schubert class  $\sigma_w$.  We will need the analogous formula on a 
general  $G/P$.  Although it is not hard to deduce such a formula from 
Chevalley's, by means of the projection  $G/B \to  G/P$,  we include a proof, 
which combines Chevalley's geometric ideas with our calculations here.  

Recall that for a simple root  $\beta$,  and a positive root  $\alpha$,  
$h_{\alpha}(\omega_{\beta}) =  n_{\alpha \beta}  
\frac{(\beta,\beta)}{(\alpha,\alpha)}$,  where  $ n_{\alpha \beta}$  is 
the coefficient of  $\alpha$  in its expansion as a positive linear combination 
of simple roots (see \S3).

\begin{lemma} [Chevalley's formula] \label{lem15}  Let  $\beta$  be in  
$\Delta \ssm \Delta_P$,  let  $u$  be in  $W/W_P$,  and let  
$\tilde{u}$  be the minimal length representative of  $u$  in  $W$.  Then 
$$
\sigma_{s_{\beta}}  \cdot   \sigma_u  =  \sum 
h_{\alpha}(\omega_{\beta})  \, \sigma_{\tilde{u} s_{\alpha}},
$$
the sum over all positive roots  $\alpha$  such that  $\ell([\tilde{u} 
s_{\alpha}]) = \ell(u) + 1$.
\end{lemma}

\begin{proof}  We must prove that, for  $v$  in  $W/W_P$  with  $\ell(v) = 
\ell(u) + 1$, $\int_X\sigma_u \cdot \sigma_{v^{\vee}}\cdot
\sigma_{s_{\beta}} = h_{\alpha}(\omega_{\beta})$ if $\tilde{u}\cdot
s_{\alpha}$ is a representative for $v$; and that $\int_X\sigma_u
\cdot \sigma_{v^{\vee}}\cdot \sigma_{s_{\beta}} = 0$ if $u$ and $v$
are not adjacent.  Note that if $\sigma_u\cdot\sigma_v^{\vee}$ is not
zero, then $Y(u) \cap X(v)$ is not empty.  This locus is fixed by the
torus $T$, and its fixed points consist of those $x(w)$ with $u
\preceq w \preceq v$ (Lemma \ref{lem9}).  Since $\ell(u) = \ell(v) -
1$, the only fixed points are $x(u)$ and $x(v)$.  It follows that,
set-theoretically at least, this intersection is the curve $C(u,v)$.
We claim that $Y(u)$ intersects $X(v)$ properly, with multiplicity
$1$, in the curve $C(u,v)$.  In characteristic zero this follows from
Lemma \ref{lem13}.\footnote{In arbitrary characteristic, Chevalley
argues as follows.  Since $Y(u)$ meets $X(u)$ transversally at the
point $x(u)$, and $x(u)$ is a nonsingular point on $X(v)$ (since
$X(v)$ is nonsingular in codimension $1$ and
$\operatorname{codim}(X(u),X(v)) = 1$), it follows that $Y(u)$ meets
$X(v)$ transversally at $x(u)$.}  Therefore
$$
\sigma_u\cdot\sigma_{v^{\vee}} = [Y(u)] \cdot [X(v)] = [C(u,v)] = d(\alpha) = 
\sum_{\beta}\, h_{\alpha}(\omega_{\beta}) \, \sigma(s_{\beta}),
$$
the last by formulas in \S3.
\end{proof}
 
\begin{lemma} \label{lem16}  If  $u \preceq v$  in  $W/W_P$,  there is a  $w$  
in  $W/W_P$  such that  $\sigma_u \cdot \sigma_w$  contains  $\sigma_v$  
with positive coefficient.
\end{lemma}

\begin{proof}  This is trivial if  $u = v$.  If  $\ell(u) = \ell(v) - 1$,  and  
$\tilde{u}$  is the minimal length representative of  $u$,  there is an  
$\alpha$  in  $R^{+}$  such that  $\tilde{u}\cdot s_{\alpha}$  is a 
representative of  $v$  (see \cite{hu:rg},\S5.11),  and  $\alpha$  is not in  
$R_P^{+}$   since  $u \neq   v$.  In this case, if we choose  $\beta$  so that  
$h_{\alpha}(\omega_{\beta}) \neq   0$,  then  $\sigma_{s_{\beta}}  \cdot  
\sigma_u$  contains  $\sigma_v$  with positive coefficient by Lemma 
\ref{lem15}.  In general, induct on  $\ell(v) - \ell(u)$,  by choosing  $u'$  
not equal to  $u$  or  $v$  with  $u \prec u' \prec v$.  If  $\sigma_{w(1)} 
\cdot \sigma_u$  contains  $\sigma_{u'}$ with positive coefficient, and  
$\sigma_{w(2)}\cdot\sigma_{u'}$  contains  $\sigma_v$  with positive 
coefficient, some  $\sigma_w$  that appears in  
$\sigma_{w(1)}\cdot\sigma_{w(2)}$  must have  $\sigma_v$  occurring in  
$\sigma_w\cdot\sigma_u$  with positive coefficient.
\end{proof}

\section{The Theorem}

The classes  $q^d\sigma_w$,  as  $d$  varies over degrees, and  $w$  varies 
over  $W/W_P$,  form a basis for the quantum cohomology ring  
$QH^{\star}(X)$  over  ${\mathbb{Z}}$.  Given any element  $\tau$  in  
$QH^{\star}(X)$,  we say that  $q^d$  {\bf occurs} in  $\tau$  if the 
coefficient of  $q^d\sigma_w$  is not zero for some  $w$.  When  $\tau$  is 
a product of Schubert classes, we know that all such coefficients are 
nonnegative.  For example,  $q^d$  occurs in  $\sigma_u \star \sigma_v$  
exactly when there is a  $w^{\vee}$  for which the GW-invariant  $\langle 
\sigma_u,\sigma_v,\sigma_{w^{\vee}}\rangle_d$  is positive.   

We now come to our main result.  The theorem gives three equivalent 
criteria for a degree to be minimal, while the proof shows these are 
equivalent to eight other related criteria.  

\begin{theorem} Let  $u$  and  $v$  be in  $W/W_P$,  and let  $d$  be a 
degree.  The following are equivalent:
\begin{enumerate}
\item[(1)]  There is a degree  $c \leq d$  such that  $q^c$  occurs in  
$\sigma_u \star \sigma_v$.

\item[(2)]  There is a chain of degree  $c \leq d$  between  $u$  and  $v$.

\item[(3)]  There is a morphism  $\varphi : {\mathbb{P}}^1 \to  X$  with  
$\varphi_{\star}[{\mathbb{P}}^1] \leq d$  such that 
$\varphi({\mathbb{P}}^1)$  meets  $Y(u)$  and  $X(v^{\vee})$.  
\end{enumerate}
\end{theorem}

\begin{proof}  We first state the eight equivalent conditions, and then we 
construct enough implications to show that each of the eleven implies the 
others.
\begin{enumerate}
\item[(4)]  There is a degree $c \leq d$,  a  $u' \succeq u$,  and a  $v' \succeq 
v$  such that  $q^c$  occurs in  $\sigma_{u'} \star \sigma_{v'}$.

\item[(5)]  There is a  $\tau$  in  $QH^{\star}(X)$  such that  $q^d$  occurs 
in  $\sigma_u \star \sigma_v \star \tau$.

\item[(6)]  The same as in (5), but with  $\tau = \sigma_{w_1} \star \ldots 
\star \sigma_{w_r}$,  for some  $w_1, \ldots , w_r$  in  $W/W_P$.  

\item[(7)]  There is a sequence  $C_0, \ldots , C_r$  of $T$-invariant curves 
on  $X$,  with  $C_0$  meeting  $Y(u)$  and  $C_r$  meeting  
$X(v^{\vee})$,  with  $C_{i-1}$  meeting  $C_i$  for  $1 \leq i \leq r$,  and 
with  $\sum_{i = 0}^r [C_i] \leq d$.

\item[(8)]  There is a connected curve  $C$  in  $X$  with  $[C] \leq d$,  
meeting  $Y(u)$  and  $X(v^{\vee})$.

\item[(9)]  There is an  $n\geq 3$  and  $w_3, \ldots , w_n$  in  $W/W_P$,  
and a  $c \leq d$,  such that, with  $e_i  : \overline{M}_{0,n}(X,c) \to  X$  
the evaluation maps, and  
$f  : \overline{M}_{0,n}(X,c) \to  \overline{M}_{0,n}$  the forgetful map,
$$
f_{\star} \bigl( {e_1}^{\star}(\sigma_u) \cdot  {e_2}^{\star} (\sigma_v) 
\cdot {e_3}^{\star}(\sigma_{w_3}) \cdot \ldots \ {e_n}^{\star} 
(\sigma_{w_n})\bigr) \, =\,  \kappa \cdot 1
$$
in  $H^0(\overline{M}_{0,n})$,  with  $\kappa > 0$.

\item[(10)]  There is a  $w$  in  $W/W_P$  and a  $c \leq d$  such that
$$
\int_{\overline{M}_{0,3}(X,c)} \, {e_1}^{\star} (\sigma_u) \cdot 
{e_2}^{\star} (\sigma_v) \cdot {e_3}^{\star} (\sigma_w)   \neq    0.
$$
\item[(11)]  There is a  $c \leq d$  such that the locus  $E(u,v^{\vee}) = 
{e_1}^{-1}(Y(u))  \cap  {e_2}^{-1}(X(v^{\vee}))$ is not empty in 
$\overline{M}_{0,3}(X,c)$.
\end{enumerate}

\vspace{12pt}

(1) $\Leftrightarrow$ (4).  The implication (1) $\Rightarrow$ (4) is trivial, 
so assume (4).  Since  $u \preceq u'$,  it follows from Lemma \ref{lem16} that 
there is a  $u''$  so that  $\sigma_{u'}$  occurs in the classical product  
$\sigma_u\cdot\sigma_{u''}$  with positive coefficient.  Take similarly  $v''$ 
for  $v \preceq v'$.  The fact that all coefficients of all quantum products of 
Schubert classes are nonnegative implies that $\sigma_u \star \sigma_{u''}  
\star \sigma_v \star \sigma_{v''}$  contains all the terms that occur in  
$\sigma_{u'} \star \sigma_{v'}$,  so it must contain some  $q^c\cdot \tau$,  
$\tau \neq   0$.  But this is a product of  $\sigma_u \star \sigma_v$  and a 
nonnegative combination of powers of    $q$'s   times Schubert classes, so  
$\sigma_u \star \sigma_v$  must contain some  $q^e$,  with  $e \leq c$.

(6) $\Rightarrow$ (5) and (5) $\Rightarrow$ (4) and (10) $\Rightarrow$ 
(9) are trivial. 

(1) $\Rightarrow$ (10).  If  $\sigma_u\star\sigma_v$  contains  
$q^c\sigma_w$  with positive coefficient  $\kappa$,  then  $f_{\star}\bigl( 
{e_1}^{\star}(\sigma_u) \cdot {e_2}^{\star}(\sigma_v) \cdot 
{e_3}^{\star}(\sigma_{w^{\vee}}) \bigr)  = \kappa\cdot 1$  in  
$H^0(\overline{M}_{0,3})$,  and  $\overline{M}_{0,3}$  is a point.

(9) $\Rightarrow$ (6).  (9) implies that  $q^c\sigma_{{w_n}^{\!\smaller\vee}}$  
occurs with coefficient  $\kappa$  in  
$\sigma_u\star\sigma_v\star\sigma_{w_3} \star \ldots \star 
\sigma_{w_{n-1}}$.

(7) $\Rightarrow$ (2).  This follows from the correspondence between 
$T$-invariant curves and pairs of adjacent elements of  $W/W_P$  (\S4).  

(2) $\Rightarrow$ (3).  A chain of degree  $c$  between  $u$  and  $v$  
corresponds to a chain of $T$-invariant curves between  $x(u)$  and  $x(v)$.  
We may assume that no curve appears more than once, since removing 
duplicates only decreases the degree.  This corresponds to a point  $\zeta = 
(C, p_1, p_2,\psi)$  in  $\overline{M}_{0,2}(X,c)$,  with $\psi   :  C \to  
X$  an embedding, $\psi(p_1) = x(u)$,  and $\psi(p_2) = x(v)$.  If  $c  = 0$,  
we are in the classical case, and $ u \preceq v^{\vee}$,  so  $x(u)$  is in  $Y(u) 
\cap X(v^{\vee})$,  and the constant map from  ${\mathbb{P}}^1$  to  $x(u)$  
satisfies the conditions of (3).  We may therefore assume that  $c > 0$.  By 
Lemma \ref{lem13},  $E(u,v^{\vee})$  is not contained in any boundary 
component.  It therefore contains a point  $({\mathbb{P}}^1,  p_1, p_2, 
\varphi)$,  and thus we have a map  $\varphi : {\mathbb{P}}^1 \to  X$  with  
$\varphi_{\star}[{\mathbb{P}}^1] = c$  and  $\varphi(p_1)$  in  $Y(u)$  and  
$\varphi(p_2)$  in  $X(v^{\vee})$.   

(3) $\Rightarrow$ (8).  With  $\varphi : {\mathbb{P}}^1 \to  X$  as in (3), 
the image curve  $\varphi({\mathbb{P}}^1)$  is connected and joins the two 
Schubert varieties, and its degree is at most  
$\varphi_{\star}[{\mathbb{P}}^1]$.

(8) $\Rightarrow$ (7).  This is a consequence of the action of the torus $T$   
on a Hilbert scheme (or Chow variety) of curves on  $X$  that join the two 
($T$-invariant) Schubert varieties, see \cite{hi:ch}.  There must be a curve in 
such a space that is fixed by  $T$,  and, being a limit of connected curves, it 
is connected.  

(10) $\Rightarrow$ (11) follows from the fact that loci representing 
cohomology classes with nonzero product must intersect, and 
${e_1}^{\star}(\sigma_u) \cdot {e_2}^{\star}(\sigma_v)$   lives on the locus  
$E(u,v^{\vee})$.

(11) $\Rightarrow$ (7).  Since the torus preserves  $E(u,v^{\vee})$,  there 
must be a point  $\zeta$  in  $E(u,v^{\vee})$  that is fixed by  $T$.  Writing  
$\zeta = (C, p_1, p_2, p_3, \varphi)$,  the image of each irreducible 
component of  $C$  must be a $T$-invariant curve in  $X$,  and the image of 
each  $\varphi(p_i)$  must be a point fixed by  $T$;  and we must have  
$\varphi(p_1)$  in  $Y(u)$  and  $\varphi(p_2)$  in  $X(v^{\vee})$.  Since 
the image is connected one can extract from it a chain of $T$-invariant 
curves from  $\varphi(p_1)$  to  $\varphi(p_2)$,  and the degree of this 
chain is at most  $\varphi_{\star}[C] = c$.  

(2) $\Rightarrow$ (9).  Suppose we have a chain of degree  $c \leq d$.  
Discarding extra curves in the chain, we may assume it is minimal.  Take  
$n$  larger than  $\sum n_{\beta} c_{\beta} + 2$.  The locus  
$E(u,v^{\vee})$  in  $\overline{M}_{0,n}(X,c)$  has pure codimension  
$\ell(u) + \ell(v)$,  and it meets the open set  $M_{0,n}(X,c)$,  by Lemma 
\ref{lem13}.  Set  
$$
 e' = e_3 \times \ldots \times e_n \,  :  \, \overline{M}_{0,n}(X,c)  \to   
X^{n-2}.
$$
Take a point  $p$  in  $M_{0,n}$,  i.e., choose  $n$  distinct points in  
${\mathbb{P}}^1$.  By Lemma \ref{lem14}, the restriction of  $e'$  to  
$f^{-1}(p) \cap M_{0,n}(X,c)$  is a finite to one mapping.  It follows that  
$(e')_{\star}[f^{-1}(p) \cap E(u,v^{\vee})] \neq   0$.  Hence there are  $w_3, 
\ldots , w_n$  in  $W/W_P$  such that  
$$
(e')_{\star}[f^{-1}(p) \cap E(u,v^{\vee})]  \,  \cdot  \, \bigl( \sigma_{w_3}  
\times \ldots \times \sigma_{w_n} \bigr)  =  \kappa\cdot[\operatorname{point}],
$$
for some  $\kappa \neq   0$.  This means that  
$$	   f_{\star} \bigl( {e_1}^{\star}(\sigma_u) \cdot 
{e_2}^{\star}(\sigma_v) \cdot {e_3}^{\star}(\sigma_{w_3}) \cdot \ldots 
\cdot {e_n}^{\star}(\sigma_{w_n}) \bigr)  =  \kappa \cdot 1  
$$
in  $H^0(\overline{M}_{0,n})$,  as required.  		
\end{proof}

\section{Peterson's quantum Chevalley formula in the parabolic case}

The quantum Chevalley formula in $G/P$ gives the formula for a quantum
product $\sigma_{s_{\beta}} \star \sigma_u$, for $\beta$ in $\Delta
\ssm \Delta_P$ and $u$ in $W/W_P$.  It starts with the classical
product $\sigma_{s_{\beta}} \cdot \sigma_u$ given in Lemma
\ref{lem15}.  The terms with $q^d$ for positive degrees $d$ have a
similar combinatorial description.  This had been proved for the
Grassmannian $\operatorname{Gr}(r,n)$ \cite{be:qs} and the complete
flag manifold $\operatorname{Fl}(\mathbb{C}^n)$ \cite{fgp:qs}.  For a
positive root $\alpha$, we use the notation $n_{\alpha}$ for
$\int_{C_{\alpha}} c_1(T_X)$ as in Lemma \ref{lem6}; and, for a simple
root $\beta$, $h_{\alpha}(\omega_{\beta})$ as before Lemma \ref{lem2}.
For $P = B$, this formula was stated by Peterson \cite{pe:mi}.

\begin{theorem} [Quantum Chevalley Formula] \label{qchev} For  $\beta$  in  $\Delta 
\ssm \Delta_P$,  $u$  in  $W/W_P$,  with  $\tilde{u}$  its 
minimal length representative in  $W$,
$$
\sigma_{s_{\beta}} \star \sigma_u  =  \sum_{\alpha} 
h_{\alpha}(\omega_{\beta})  \sigma_{\tilde{u} s_{\alpha}} +  
\sum_{\alpha} q^{d(\alpha)} h_{\alpha}(\omega_{\beta})  \sigma_{\tilde{u} 
s_{\alpha}},
$$
the first sum over roots $\alpha$ in $R^{+} \ssm R_P^{+}$ for which
$\ell(v) = \ell(u) + 1$, where $v$ is the coset of $\ti{u} s_\alpha$
in $W/W_P$, and the second sum over roots $\alpha$ in $R^{+} \ssm
R_P^{+}$ for which $\ell(v) = \ell(u) + 1 - n_{\alpha}$.
\end{theorem}

\begin{proof}  Let  $q^d\sigma_v$  be a term appearing in  
$\sigma_{s_{\beta}} \star \sigma_u$ with nonzero coefficient $\kappa$,
i.e.,
$$
\kappa \, = \,  \langle  \sigma_u, \sigma_{v^{\vee}}, 
\sigma_{s_{\beta}}\rangle_d  
= \int_{\overline{M}_{0,3}(X,d)} {e_1}^{\star}(\sigma_u)\cdot 
{e_2}^{\star}(\sigma_{v^{\vee}})\cdot {e_3}^{\star}(\sigma_{s_{\beta}}) \,  
> \, 0.  
$$
It suffices to consider the case $d \neq 0$, since the classical case
is covered by Chevalley's formula \ref{lem15}.  Set $E = E(u,v)$ in
$\overline{M}_{0,3}(X,d)$.  By the transversality Lemma \ref{lem13},
$E$ is reduced, locally irreducible, purely 1-dimensional, with each
component meeting $M_{0,3}(X,d)$, and with
$$
[E]  =  {e_1}^{\star}([Y(u)]) \cdot  {e_2}^{\star}([X(v)])  =  
{e_1}^{\star}(\sigma_u)\cdot {e_2}^{\star}(\sigma_{v^{\vee}}).
$$
So  $\kappa$  is the coefficient of  $\sigma(s_{\beta})$  in the class  
$(e_3)_{\star}[E]$.   

We claim that $e_1(E) = x(u)$ and $e_2(E) = x(v)$.  Indeed, each of
the loci $e_i(E)$ is $T$-invariant, and if $e_1(E) \neq x(u)$, then
$e_1(E)$ would contain some other $T$-fixed point $x(u')$, with $u
\prec u'$, $u \neq u'$.  But then $E$ would be contained in
$E(u',v)$, which (by transversality again) has larger codimension than
$E$, a contradiction; and similarly if $e_2(E)$ contains $x(v')$, with
$v' \prec v$.

Note that  $\operatorname{codim}( E) = \ell(u) + \operatorname{dim}( X) - 
\ell(v) = \operatorname{dim}(\overline{M}_{0,3}(X,d)) - 1 = 
\operatorname{dim}( X) + \int_d  c_1(T_X)$,  so  
$$
\ell(v)  =  \ell(u) + 1 - \int_d \, c_1(T_X) .
$$
In particular, $\ell(v) < \ell(u)$, since, by
Lemma \ref{lem6}, $\int_d c_1(T_X) \geq 2$.

A point of  $E$  not in the boundary has the form  $\zeta = ({\mathbb{P}}^1, 
p_1, p_2, p_3, \varphi)$,  with  $\varphi : {\mathbb{P}}^1 \to  X$,  
$\varphi_{\star}[{\mathbb{P}}^1] = d$,  and  $\varphi(p_1) = x(u)$  and  
$\varphi(p_2) = x(v)$.   Given such a map  $\varphi$,  one can produce a 
curve in  $E$  containing this point by varying where  $\varphi$  maps  
$p_3$  in the curve  $\varphi({\mathbb{P}}^1)$,  i.e., varying  $\varphi$  by  
$\varphi\circ\vartheta$,  where  $\vartheta$  is an automorphism of  
${\mathbb{P}}^1$  that fixes  $p_1$  and  $p_2$.  It follows that each 
irreducible component  $Z$  of  $E$  must consist generically of such maps.  
It follows from this that  $\varphi({\mathbb{P}}^1) = e_3(Z)$.  Since  
$e_3(Z)$  is invariant by  $T$,  it follows that  $\varphi({\mathbb{P}}^1)$  
is  $T$-invariant, and, since it contains  $x(u)$  and  $x(v)$,  we must have  
$\varphi({\mathbb{P}}^1) = C(u,v)$.  In particular,  $u$  and  $v$  must be 
adjacent. 

We claim next that  $\varphi$  maps  ${\mathbb{P}}^1$  isomorphically 
onto  $C(u,v)$.  If not,  $\varphi_{\star}[{\mathbb{P}}^1] = k[C(u,v)]$,  for  
$k > 1$.  Consider the corresponding locus  $E' = E(u,v)$  in  
$$\overline{M}' = \overline{M}_{0,3}(C(u,v),k[C(u,v)])$$  
consisting of maps to $C(u,v)$ of degree $k$
that take  $p_1$  to  $x(u)$  and  $p_2$  to  $x(v)$.  The codimension of  
$E'$  in  $\overline{M}'$  is at most  $2$  (in fact, equal to $2$ by Lemma 
\ref{lem13} applied to the variety  $C(u,v)$  in place of  $X$).  But the 
dimension of  $\overline{M}'$  is  $2k+1$,  and since  $E' \subset E$,  we 
must have  $2k-1 \leq 1$,  i.e.,  $k = 1$.  

It follows that  $d = [C(u,v)] = d(u,v) = d(\alpha)$,  where  $\tilde{u} 
s_{\alpha}$  is a representative of  $v$.  It also follows that each component  
$Z$  of  $E$  is equal to the locus  $E'$  described in the preceding 
paragraph.  In particular, there is only one irreducible component  $E = E'$  
of  $E$. (This  $E$  can be realized from the blow-up of  $C(u,v)  \times  
C(u,v)$  along the two points  $x(u)  \times  x(u)$  and  $x(v)  \times  x(v)$;  
the exceptional divisors become the extra factors in curves corresponding to 
the two boundary points of  $E$.)  It also follows from this description that  
$e_3$  maps  $E$  isomorphically onto  $C(u,v)$.  Therefore  
$(e_3)_{\star}[E] = [C(u,v)] = d(\alpha)$,  so  $\kappa$  is the coefficient of  
$\sigma(s_{\beta})$  in  $d(\alpha)$,  which is  
$h_{\alpha}(\omega_{\beta})$,  as we saw in Lemma \ref{lem5}.  This 
argument, run backwards, shows that each such  $v$  does occur, and 
completes the proof of the theorem.
\end{proof}

For the classical flag manifold  $\operatorname{Fl}(\mathbb{C}^n)$,  
the positive roots and 
reflections correspond to  $\alpha = (a,b)$,  $1 \leq a < b \leq n$,  and  
$h_{(a,b)}(\omega_{(r,r+1)})$  is  $1$  if  $a \leq r < b$  and is  $0$  
otherwise.  So one recovers the quantum Monk formula of \cite{fgp:qs}.

Peterson's generalization of the Chevalley formula allows the quantum
products of classes in the subalgebra generated (using the quantum
product) by $H^2(G/P,\Z) \otimes \Z[q]$ to be computed recursively.
S. Fomin has pointed out that this subalgebra is the full quantum 
cohomology ring $QH^*(G/P)$ whenever the classical ring $H^*(G/P)$ 
is generated by $H^2(G/P)$.  (The proof is by induction on the degree.) 
This holds on the flag varieties $G/B$.  

Peterson has also given an explicit formula for any Gromov-Witten
invariant on any $G/P$ as another Gromov-Witten invariant on the
corresponding $G/B$.  A proof of this, which uses some of the ideas of
this paper, is given in \cite{wo:pe}.



\vskip .2in


\end{document}